\documentclass[final,pdftex]{siamltex} 
\pdfoutput=1

\usepackage{amssymb,amsmath,amsfonts}
\usepackage{graphicx}
\usepackage{exscale}
\usepackage{rotating}
\usepackage{tikz}
\usetikzlibrary{arrows,snakes,backgrounds,calc}

\usepackage{color}

\usepackage{bm}
\usepackage{time}

\usepackage[normalem]{ulem}
\usepackage[nooneline,hang]{subfigure}

\usepackage{grffile}

\newcommand\Fig[1]{Fig.~\ref{fig:#1}}

\newcommand\Eq[1]{Eq.~(\ref{eq:#1})}

\newcommand{\ie}{{i.\,e., }}
\newcommand{\eg}{{e.\,g., }}
\newcommand{\cf}{{cf.\ }}

\def\6{\partial}

\def\Re{\ensuremath{\textrm{Re}}}

\title{Unidirectionally Coupled Map Lattices with Non-Linear Coupling:
  Unbinding Transitions and Super-Long Transients}
\author{Christian Marschler%
  \thanks{Technical University of Denmark, Department of
    Applied Mathematics and Computer Science, Matematiktorvet 303B,
    DK-2800 Kgs. Lyngby, Denmark, and Max Planck Institute for Dynamics
    and Self-Organization, Fassberg~17, D-37077 G\"ottingen,
    Germany 
    ({\tt c.marschler@mat.dtu.dk}).}%
  \and J\"urgen Vollmer%
  \thanks{Max Planck Institute for Dynamics and Self-Organization,
    Fassberg~17, D-37077 G\"ottingen, Germany, and Institute for
    Nonlinear Dynamics, Faculty of Physics, Georg August University,
    D-37077 G\"ottingen, Germany ({\tt juergen.vollmer@ds.mpg.de}).}
}
\begin{document}
\maketitle

\newcommand{\slugmaster}{}

\begin{abstract}
  Recently, highly resolved experiments and simulations have provided
  detailed insight into the dynamics of turbulent pipe flow. 
  This has revived the interest to identify mechanisms that generate
  chaotic transients with super-exponential growth of lifetime as a
  function of a control parameter, the Reynolds number for pipe flow,
  and with transitions from bounded chaotic patches to an invasion of
  space of irregular motion.  Dynamical systems models are unique
  tools in this respect because they can provide insight into the
  origin of the very long life time of puffs, and the dynamical mechanism
  leading to the transition from puffs to slugs in pipe flow.
  The present paper contributes to this enterprise by introducing a
  unidirectionally coupled map lattice. It mimics three of the
    salient features of pipe-flow turbulence: (i)~the transition from laminar flow to
  puffs, (ii)~a super-exponential scaling of puff
  lifetime, and (iii)~the
  transition from puffs to slugs by an unbinding transition in an
  intermittency scenario. In our model all
  transitions and scalings are theoretically described from a dynamical
  systems point of view.
 \end{abstract}



\section{Introduction}

More than 100 years ago Reynolds \cite{Reynolds1883} established that
the qualitative behavior of pipe flow is described by a single control
parameter, the Reynolds number \Re, while details of trajectories
depend on the specific initial conditions.
For all $\Re$ a parabolic flow profile appears to be stable
\cite{Meseguer2003,WygnanskiChampagne1973}. However, for $1750
\lesssim \Re \lesssim 2200$ there emerges another solution
\cite{MellibovskyMeseguerSchneiderEckhardt2009}: long-living
chaotic transients of a fixed finite-size, \emph{turbulent puffs}, that
are advected down the pipe with the mean flow.  Their dynamics
shows sensitive dependence on initial conditions
\cite{FaisstEckhardt2004}, and for any given Reynolds number the
lifetime of individual puffs is exponentially distributed
\cite{FaisstEckhardt2004,HofWesterweelSchneiderEckhardt2006,%
  PeixinhoMullin2006,
  HofLozarKuikWesterweel2008}.
There appears to be no critical Reynolds number where the lifetime of
the puffs diverges. Rather, recent experimental data suggest
\cite{HofLozarKuikWesterweel2008} that the mean lifetime $\tau$ grows
super-exponentially like $\log\log\tau \sim \Re$, in accordance with
direct numerical simulations \cite{Avila2010}.
For $\Re \gtrsim 2200$ there is an increasing probability for
puff splitting
\cite{WygnanskiChampagne1973,MoxeyBarkley2010,Avila2011}.  A critical
point for this transition has been shown to emerge at $\Re = 2040$
\cite{Avila2011}, marking the cross-over between two processes: decay
of turbulence by the decay of puffs, and its spreading by puff
splitting.  Eventually the splitting leads to a constant growth of
the turbulent regions in the pipe, \emph{turbulent
slugs}. Numerical simulations for pipe flow reveal highly non-linear,
\ie non-diffusive, coupling in the downstream
direction~\cite[Fig.~18]{DuguetWillisKerswell2010}, 
where intermittently very strong azimuthal vorticity perturbations
are generated close to the pipe wall. Subsequently, they propagate to 
the center of the pipe where they are accelerated and enhance 
turbulence downstream the pipe.

The important role played by the nonlinear forward coupling in pipe flow 
motivates us to study the consequences of non-linear
forward coupling from a dynamical systems perspective.
We establish a coupled map lattice (CML) where puffs are instances of
super-long transients in dynamical systems \cite{TelLai2008}.
CMLs are dynamical systems where the dynamic quantities take
continuous values defined on a discrete set of lattice sites, and
where these values are updated at discrete times according to a
deterministic rule.
Our CML features a
transiently
chaotic on-site dynamics and a non-linear, unidirectional coupling
between lattice sites. 

CMLs have been studied extensively for diffusively-coupled, chaotic
on-site dynamics \cite{Bunimovich1992,Bunimovich1988}. In particular,
the study of spatio-temporal intermittency, chaos and supertransients
\cite{kaneko1986,Kaneko1989,Kaneko1990,ChateManneville1988,GrassbergerSchreiber1991,Willeboordse1994}
and the control of chaos \cite{Gang1994,Parmananda1997,Parekh1998}
gained much attention. Applications reach from front motion in
hydrodynamic models \cite{Pomeau1986}, to modelling of traffic flow
\cite{Yukawa1995} and evolution of genetic sequences \cite{Cocho1991}.
In Refs.~\cite{Rudzick1996,Rudzick1997} advectively-diffusively
coupled map lattices have been studied as a
model for open flow, addressing in particular the relation to
turbulent pipe flow and to spatio-temporal intermittency
\cite{PomeauManneville1980,kaneko1986,ChateManneville1988,Pomeau1989}.
A recent model along this line picked up this approach to gain insight
into the transition from puffs to slugs in turbulent pipe flow
\cite{Barkley2011}.

Unidirectionally
  coupled map lattices (UCML) were introduced in Ref.~\cite{Rudzick1996} in order to address the
transition to turbulence in open flows.  UCMLs provide a minimalistic
approach to understand the dynamics and transitions, \ie bifurcations
in the dynamics of turbulent pipe flow.  However, in contrast to
the models considered previously \cite{Rudzick1996,Barkley2011} our
UCML has a non-linear coupling. As a consequence the model is
complicated enough to show the complex behavior observed in turbulent
pipe flow. On the other hand, it is simple enough so that we can
analytically discuss the bifurcation points for the transitions from
laminar flow to turbulent puffs (Section \ref{sec:puffs}) and from
puffs to slugs (Section \ref{ssec:laminarpuffs}).

The paper is structured as follows. In Section \ref{sec:model} we
introduce our model.  Section \ref{sec:numres} gives numerically
computed sample trajectories for different parameter values to get an
idea about the different solutions of the system.
Section~\ref{sec:results} explains the main results and compares
analytical and numerical findings. In particular, the transitions from
laminar pipe flow to puffs, and from puffs to slugs are presented in
Sections \ref{sec:puffs} and \ref{ssec:laminarpuffs}, respectively.  A
lifetime analysis of puffs is presented in
Section~\ref{ssec:lifetime}, and the implications of our model to the
qualitative understanding of pipe-flow turbulence is further discussed
in Section~\ref{sec:discussion}.  The paper concludes with a summary
and an outlook to future research directions in Section
\ref{sec:conclusion}.

\section{Model}
\label{sec:model}
In this section, we define our unidirectionally
coupled map lattice (UCML). It is defined on a one-dimensional lattice
with sites $i \in \lbrace 1,\ldots,m \rbrace$.

We represent the state of the system at site $i$ and
time $t\in \mathbb{N}$ by
variables $x_t^i \in \mathbb{R}$ which mimick the turbulent kinetic energy of the flow.
For a laminar flow, \ie the stable stationary solution for the pipe flow
problem, this variable takes values $x^i_t = 0$. 

\begin{figure}[t]
  \centering 
%
\subfigure[Sketch of the unidirectionally coupled map
  lattice.]{\label{fig:ucml}\includegraphics[width = .7\textwidth]{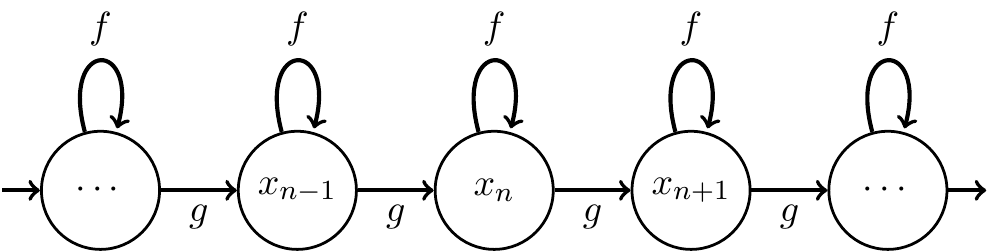}}
%
    \subfigure[On-site dynamics $f$.]{\label{fig:f}\includegraphics[width = 0.45\textwidth]{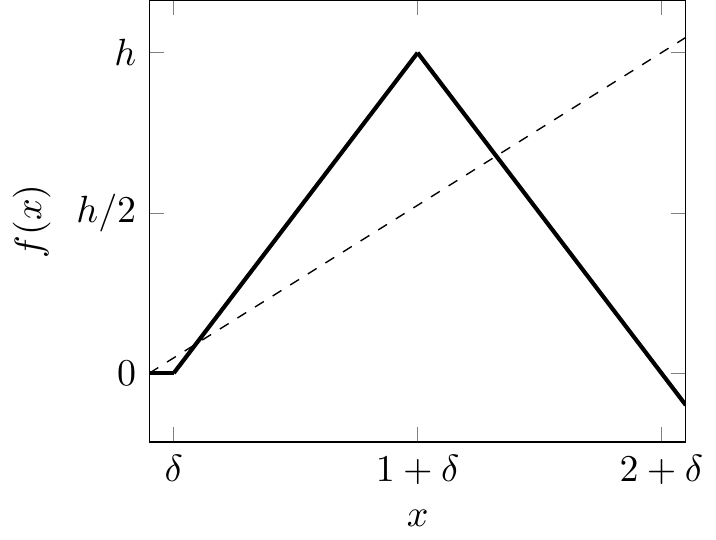}}
    \hfill 
%
    \subfigure[Coupling function $\alpha\, g$.]{\label{fig:g}\includegraphics[width = 0.45\textwidth]{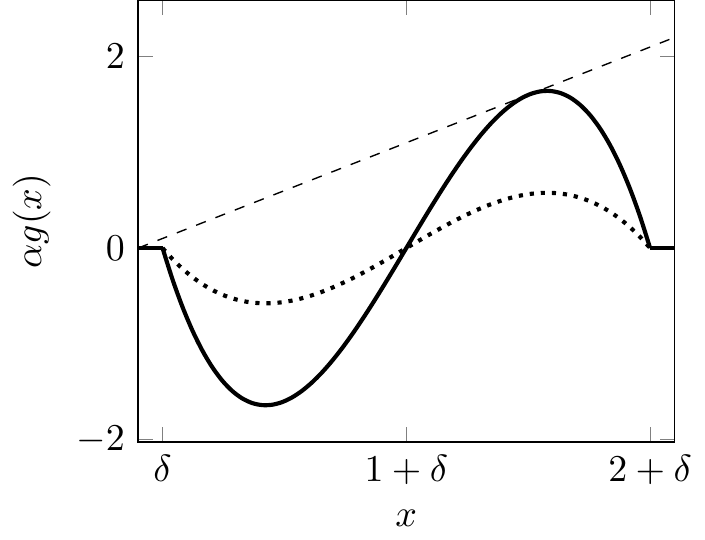}}
    \caption{\textup{(a)} Sketch of our unidirectionally coupled map
      lattice. The on-site dynamics is given by the function $f$, and
      the coupling function $g$ determines the forward coupling,
      mimicking the downstream of pipe flow. \textup{(b)} The on-site
      dynamics $f$ is a shifted tent map \Eq{f(x)}, which is cut off
      at $\delta$. Intersections with the diagonal (dashed line)
      denote fixed points. For $h>2$, almost all trajectories are
      transients to the globally stable fixed point at
      $x=0$. \textup{(c)} The coupling function $\alpha\, g$ for the
      forward coupling takes the form specified in \Eq{g(x)}. For
      $\alpha < \alpha_{sn}$ (dotted curve), it has a single fixed
      point $x=0$. At $\alpha = \alpha_{sn}$ (solid curve), another
      pair of fixed points is born in a saddle-node bifurcation.}
  \label{fig:fandg}
\end{figure}

The evolution of the state variables $x_t^{i}$ follows a
coupled map lattice dynamics 

\begin{equation}
\label{eq:model}%
  x_{t+1}^i =  \alpha \: g(x_t^{i-1}) + f(x_t^i),
\end{equation}
where $f$ generates an on-site dynamics, and $\alpha\, g$
is a \emph{unidirectional} coupling of sites $i$ and $i-1$ with coupling strength
$\alpha$ (cf. \Fig{fandg}).

The substantial difference to previous
approaches \cite{Barkley2011,ChateManneville1988,GrassbergerSchreiber1991,%
  Kaneko1990, Rudzick1996,GinelliLiviPoliti2002} lies in our choice of
the coupling.
In contrast to other coupled map lattice models, which address
diffusive and advective coupling (see \eg
\cite{ChateManneville1988,Kaneko1989,Kaneko1990,Rudzick1996,Barkley2011}), we adopt
a non-linear unidirectional coupling (cf.\ \Fig{g})
\begin{equation}
  g(x) = 
  \begin{cases}
    - \frac{3}{2} \, 
    \left( x - \delta \right)
    \left( x - 1 -\delta \right)
    \left( x - 2 -\delta \right) \, , 
& \delta \leq x < 2+\delta \\[1mm]
\;\;\; 0 \, , & \text{else.}
\end{cases}
\label{eq:g(x)}
\end{equation}%
In our choice to provide only a forward coupling we 
follow Ref.~\cite{Rudzick1996,Rudzick1997}. The choice of 
a nonlinear coupling function is
motivated by the numerical results shown in Fig.~18 of \cite{DuguetWillisKerswell2010}, 
where spatial intermittency in pipe flow is connected to the fast 
forward propagation of vorticity perturbations generated at the edge of the pipe.
This finding points to a strong, nonlinear downstream coupling 
in pipe flow.
It
motivates us to choose a non-linear coupling along the lattice that
provides an intermittent dynamics close to the saddle-node
bifurcation of $\alpha\, g$ (cf.\ \Fig{g} and \Eq{sn} below).
Hence, the action of the function $g$ models the propagation of
turbulence in puffs and slugs. The specific choice of $g$ is
  generic as long as the maximum in the right domain
  $x\in[1+\delta,2+\delta]$ allows for an intermittency scenario via a
  saddle-node bifurcation.

In line with previous approaches \cite{ChateManneville1988,GinelliLiviPoliti2002,Barkley2011}
the on-site dynamics is generated by a modified tent map (\Fig{f}), 
\begin{equation}
f(x) = \begin{cases}
  \;\;\; h \, \left( x - \delta \right) \, ,
  & \delta \leq x < 1+\delta \\[2mm]
  -h\, \left( x - 2 - \delta \right) \, , 
  & 1+\delta \leq x\\[1mm]
  \;\;\; 0 \, , 
  & \text{else} \, .
\end{cases}
\label{eq:f(x)}
\end{equation}
It is augmented
with a flat piece close to the origin $[0,\delta]$ in accordance with the
requirement that $x=0$ is a globally stable fixed point.
For $h>1$, the on-site dynamics has three fixed points
\begin{equation}
  \label{eq:fpoff}
  x_0^f = 0, \qquad x_1^f = \frac{h \delta}{h -1},\qquad x_2^f = \frac{h (2+\delta)}{1+h},
\end{equation}
where $x_0^f$ is stable, and $x_1^f$, $x_2^f$ are unstable.  The
specific form of the cutoff towards the very small values of $x$ does
not affect the scaling of the transitions to be discussed in the
following, as long as $x=0$ remains a globally stable fixed point. For
the purpose of the present study the offset is kept fixed,
$\delta=0.1$.  Follow up studies that strive for a quantitative
comparison with pipe flow will have to account for the shrinking of
the basin of attraction for the laminar flow with $\Re$ (see
\cite{HofJuelMullin03}) in real pipe flow.  Further studies of system
\eqref{eq:model} could include a dependence $\delta (\alpha,h)$,
effectively modelling $\delta (\Re)$, and reducing the number of free
parameters in our UCML.
Since $Df(x_0^f) \equiv 0$, the state $x_t^i=0$ for all $i$ is a
stable fixed point of the dynamics for all parameter values. It mimics
laminar flow. Furthermore, a turbulent state is identified with values
$x_t^i>0$. For all numerical simulations periodic boundary conditions
are enforced. During simulations, the lattice is always chosen long
enough to ensure, that the leading edge never meets the trailing edge,
\ie an infinitely-long pipe is modelled.

For $h>h_c=2$ the map $f$ shows transient chaotic behavior with an
exponential distribution of lifetime $\tau$ whose average 
\begin{equation}
  \label{eq:taus}
\tau_s = \left[ \ln \frac{h}{h_c} \right]^{-1}  
\end{equation}
diverges algebraically as $h$ approaches $h_c=2$
\cite{GrebogiOttYorke1982,Ott2002}. 

The slope $h$ of the tent map and the coupling strength $\alpha$ serve
as control parameters of the dynamics.  In Sections \ref{sec:numres}
and \ref{sec:results} we analyze the behavior of this dynamical system
in response to varying the control parameters $h$ and $\alpha$.  The
connection to pipe flow is detailed in Section \ref{sec:results} where we
discuss how the parameters can be connected to the Reynolds number,
\Re. For the discussion of the dynamical systems aspects of the
transitions this approach is more transparent, and in order to see
their physical implications it is sufficient to observe that
increasing $\Re$ amount to an increasing $\alpha$ and decreasing $h$.
We also refrain from identifying a specific frame of reference for the
velocities because it will depend in detail of the connection of the
nonlinear coupling function, $g$, and the intermittent forward
coupling in pipe flow which we do not explore with our generic choice
of $g$.

\section{Numerical Results}
\label{sec:numres}
\begin{figure}[t]
\centering
\subfigure[Turbulent puff \newline $\alpha=0.5$, $h=2.05$.]%
{\label{fig:puff}\includegraphics[width=.32\textwidth]{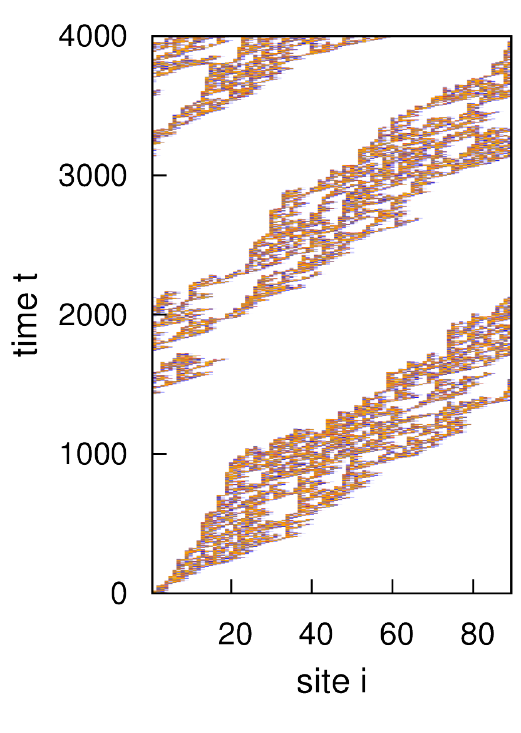}}
\hfill
\subfigure[Turbulent puff\newline $\alpha=0.8$, $h=2.10$.]%
{\label{fig:puffshort}\includegraphics[width=.23\textwidth]{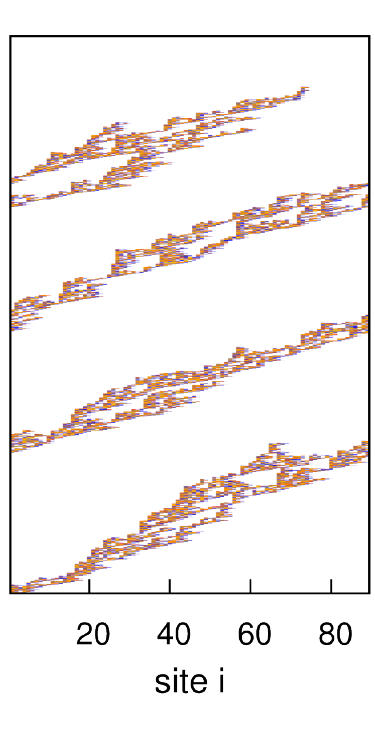}}
\hfill
\subfigure[Turbulent slug\newline $\alpha=2.8$, $h=2.10$.]%
{\label{fig:slug}\includegraphics[width=.296\textwidth]{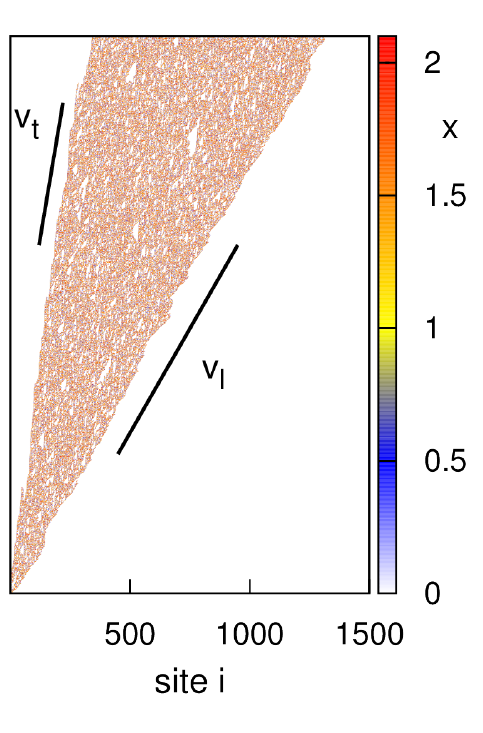}}

\caption{Space-time plots of the UCML, \Eq{model}. 
  The color represents the values of $x$ as specified by the
  color bar to the right.  
  \textup{(a)} For intermediate values of the coupling, $\alpha=0.5$, $h=2.05$,
  there are chaotic regions of a fixed
  average size that propagate with a finite speed. 
  \textup{(b)} For a different set of parameters, $\alpha = 0.8$, $h=2.10$, the lifetime of turbulent puffs can be much shorter.
  \textup{(c)} For higher coupling, $\alpha=2.8$, $h=2.10$,
  the chaotic region increases linearly in size. 
  The time axes are identical for all three plots; beware however the different system size.
  \label{fig:spacetime}}
\end{figure}

In order to provide an overview about the solutions of the UCML, we
discuss trajectories for different parameter values. Changing the
coupling parameter $\alpha$ and the height of the tent map $h$, three
qualitatively different types of solutions are observed: (i)
immediately decaying structures; (ii) propagating turbulent puffs, \ie
constant size turbulent spots that are advected downstream (cf.\
\Fig{puff} and \Fig{puffshort}); (iii) turbulent slugs, \ie growing
turbulent regions (cf.\ \Fig{slug}).

%
In \Eq{taus} the scaling of the lifetime of the on-site dynamics is stated for values
$h>h_c=2$. For $h<h_c$ chaotic on-site motion does not decay. As a consequence,
the chaotic dynamics is only spreading due to the coupling
and there will be no fluctuations where the systems relaminarizes at times. Such a
dynamics is not 
of interest to model pipe flow \cite{Avila2013}. Therefore, in the
following we only investigate solutions for $h>h_c$.
The numerical observations suggest, that there are two critical
parameter values for the coupling strength: $\alpha_P$ and
$\alpha_{sn}$. 
\begin{description}
\item For $\alpha \in [0,\alpha_P]$, the coupling is not strong enough
  to kick the adjacent lattice site into a turbulent state. A finite
  single-site initial perturbation decays unless it is exactly
  initialized at an unstable periodic orbit of $f$. The probability to
  be found in a turbulent state for $t$ successive time steps decays
  exponentially with the decay rate of the on-site dynamics,
  $\tau_s^{-1}$.
\item For $\alpha \gtrsim \alpha_P$, an initial perturbation spreads
  if $2+\delta > x>1+\delta$, \ie for $x$ where $g(x) > 0$. For unidirectional
  couplings the transient dynamics of $f$ also governs the decay of
  turbulence at the rear side of the turbulent structures.  It reaches
  a certain width. Eventually, a constant number of sites are in a
  turbulent state, and this finite-size structure travels through the
  lattice in the direction given by the coupling (\Fig{puff} and
  \Fig{puffshort}). This dynamics is reminiscent of turbulent puffs
  and it is observed numerically in the UCML, that the lifetime
  depends on the parameters. In \Fig{puffshort}, a puff solution with
  a relatively small width and short lifetime is shown. After about
  3500 time steps, the puff decays to the laminar state. For a
  different set of parameters, there are broader
  structures with a longer lifetime (\Fig{puff}). Upon increasing
  the coupling strength $\alpha$ or decreasing $h$ towards $h_c$ the
  average lifetime of puffs rises, and eventually, there is a
  transition to slug-like dynamics 
  (\Fig{slug}): the velocity at
  the leading edge, $v_l$, is higher than the velocity at the trailing
  edge, $v_t$, such that the turbulent structures broaden
  while traveling through the lattice. Hence, we find solutions
  reminiscent of turbulent slugs.
\item For $\alpha>\alpha_{sn}$, the coupling function $\alpha\, g$
  attains a non-trivial fixed point that induces ballistic propagation
  of the leading edge of the slug, $v_l = 1$.
\end{description}

In Section \ref{sec:results} we provide a detailed study of the transitions. 

\section{Analytical Treatment}
\label{sec:results}
In this section, the critical couplings $\alpha_P$ and $\alpha_{sn}$ are
computed analytically, and their dependence on $h$ is discussed (see \Fig{phase_space}).
All predictions are verified by comparing to numerically obtained results. 
We first focus on the growth speed of the turbulent region in the slug
regime, \ie on the parameter dependence of the difference $v_l - v_t$ of the propagation velocity
$v_l$ of the leading edge and the one $v_t$ of the trailing edge of
the turbulent region, as shown in 
\Fig{slug}.
 Subsequently, we
address the factors governing the lifetime of puffs
(cf.\ \Fig{puff} and \Fig{puffshort}).

\subsection{Slugs}
\label{ssec:laminarpuffs}
\begin{figure}[t]
\centering
\includegraphics[width=.85\textwidth]{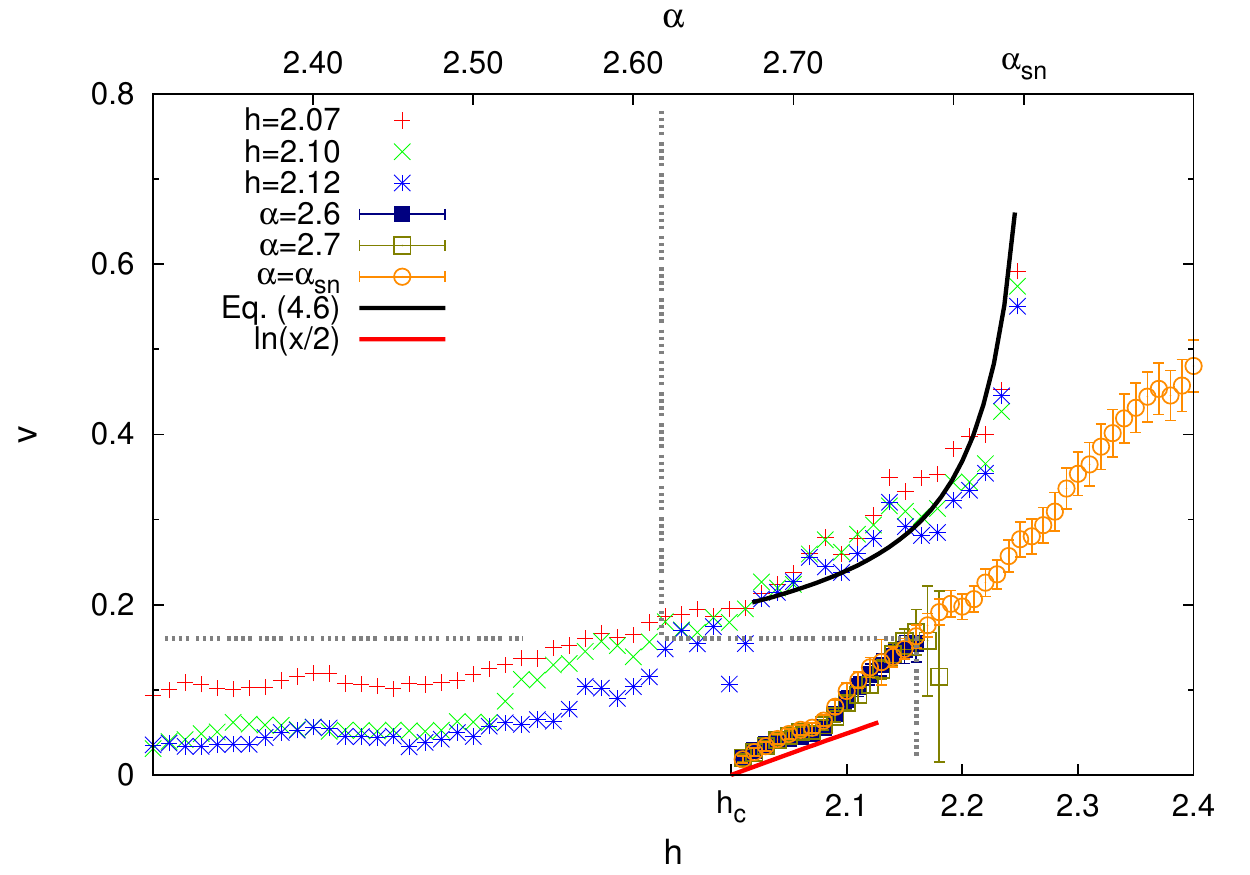}
\caption{The average propagation velocities $v_l(\alpha)$ and $v_t(h)$
  of the leading (upper left data sets, top axis) and the trailing
  edge (lower right data sets, bottom axis), respectively.
  The solid black and red lines provide analytical estimates for $v_l$
  and $v_t$, respectively, and the gray, dotted line is used to
  construct the transition line $\alpha_{PS}(h)$ from puffs to slugs
  (\cf main text for details).
  The displayed values for both velocities are obtained by averaging
  over $2\times 10^4$ initial conditions. Corresponding standard
  deviations are of the order of $0.05$ for $\alpha$ close to
  $\alpha_{sn}$, and increase to about $0.1$ for 
  $\alpha \lesssim 2.5$.
  \label{fig:slug_velocities}}
\end{figure}

Since there is only forward coupling of the dynamics, the lifetime of
the site at the trailing edge, \ie the last non-zero site of the slug
amounts to the lifetime, $\tau_s(h)$, of chaotic transients of the
on-site dynamics, $f$. The
velocity $v_t$ should hence be a function of $h$ only. This is
confirmed by the data to the lower right in \Fig{slug_velocities}.
Using \Eq{taus}, the average propagation velocity of the trailing edge
is estimated as
  \begin{equation}
    \label{eq:vtrail}
v_t \simeq \frac{1}{\tau_s(h)}  =  \ln \frac{h}{h_c}.    
  \end{equation}
  The solid red line in \Fig{slug_velocities} shows that
  \Eq{vtrail} is a lower bound for the trailing edge
    velocity. It is a sharp bound when the single-site
  lifetimes are large, \ie for $h$ close to $h_c$.
The mismatch for larger $h$ is due to jumps by several sites, when 
the right neighbors of the trailing site are already in the laminar
state. The decay of the turbulent dynamics of that site leads then to
a jump of the trailing border of the slug by more than a single site,
and hence to a higher propagation velocity. This effect becomes more
pronounced when there is a larger number of laminar sites within the
turbulent region, \ie when $\tau_s$ decreases due to an increase
of $h$.
It is evident from \Fig{spacetime} that the time scale $\tau_s$ is much 
shorter than the turbulence lifetime $\tau$. After all, turbulence is sustained 
by propagation along the system. In particular, turbulence persists when the trailing edge 
can not catch up with the leading edge.

The propagation velocity at the
leading edge of the slug, \ie the right-most site with non-vanishing $x$, is governed by
$\alpha\, g$. Hence, we expect that the velocity $v_l$ should be a
function of $\alpha$.
In order to gain insight into the propagation we start from strong coupling. For
$\alpha = \alpha_{sn} \simeq 2.844137$  
there is a saddle-node bifurcation where
a fixed point $x^f$ emerges in the coupling function $\alpha\, g$ (cf. \Fig{g}), 
\begin{equation}
\label{eq:sn}
\alpha_{sn} \ g(x^f) = x^f  \, ,
\qquad
\alpha_{sn} \ g'(x^f) = 1   \, .
\end{equation}%
Beyond this bifurcation, \ie for $\alpha > \alpha_{sn}$, the stable fixed point leads to a ballistic
propagation of the front. This behavior can be understood by analyzing
the case $\alpha = \alpha_{sn}$. If $(x_t^i,x_t^{i+1})=(x^f,0)$, then according to \Eq{model}
$x_{t+1}^{i+1} = x^f$. The turbulence propagates
to the next site in each time step; a front edge at $x=x^f$ propagates
ballistically.

Immediately before the bifurcation,
\ie for $\Delta \alpha \equiv \alpha_{sn} - \alpha \ll 1$, 
trajectories can spend a time $T$ in the vicinity of $x^f$ such that
the propagation velocity can be estimated as an average of a fast
propagation when $x \simeq x^f$ and slow propagation for other values
of $x$.

According to Little's lemma \cite{Little1961} propagation is observed
with probability 
\begin{equation}
  \label{eq:pLittle}
p(\Delta\alpha) = [1 + 1/ (\nu T)]^{-1}\, ,  
\end{equation}
where
$\nu$ is the injection rate into the vicinity of the emerging fixed
point.
\Eq{pLittle} can be understood intuitively by considering a
  2-state Markov chain, where $p$ and $q$ are the probabilities to be in the
  propagating and non-propagating state, respectively. Further, we
  denote the transition rate from state $p$ to $q$ with
  $\nu_\text{out}$ and the one from $q$ to $p$ with $\nu$. Employing
  the normalization $1 = p + q$ and the steady-state assumption $\dot
  p = \dot q = 0$, we obtain \Eq{pLittle}, where the escape rate is
  $\nu_\text{out} = 1/T$.
Moreover, the scaling of lifetimes in the intermittency scenario \cite{PomeauManneville1980,Ott2002} for the
``stickiness'' of the coordinate $x^f$, states that
\begin{equation}
\label{eq:Tint}
  T = a \; \Delta \alpha^{-\frac{1}{2}} \, ,
\end{equation}
which is the mean time to escape a narrow tunnel just before a
  saddle-node bifurcation (cf. the type-I intermittency scenario in \cite{Ott2002}).
In the present case the factor of proportionality amounts to 
$a \simeq 1.55$ and has been fitted to simulation data. $\Delta
  \alpha$ denotes the distance to the saddle-node bifurcation in $g$.
Moreover, the entrance probability to the region can be evaluated
following recent ideas on Poincar\'e recurrences
\cite{AltmannSilvaCaldas2004,Altmann2013}, 
\begin{equation}
\label{eq:nu}
\nu = \nu_c + A \exp\left(-{\Delta \alpha}/{\xi}\right)
\end{equation}
with $\nu_c=0.039$, $A=0.034$ and $\xi=0.023$ obtained by fits to
simulation data for the present map.
For $\alpha$ values below the bifurcation, 
$2.7 \lesssim \alpha < \alpha_{sn}$ and using \Eq{pLittle}, \Eq{Tint}
and \Eq{nu}, we thus obtain an explicit
expression for the spreading velocity
\begin{equation}
   v_l \simeq  p(\Delta\alpha) = \left[1 + \frac{\Delta \alpha^{\frac{1}{2}}}{a(\nu_c + A \exp\left(-{\Delta
       \alpha}/{\xi}\right)) )}\right]^{-1}.
\label{eq:frontvelocity}
\end{equation}
It is indicated by a solid black line in \Fig{slug_velocities}
and helps us to understand the leading edge velocity and thereby
  the transition from puffs to slugs. 

\begin{figure}[t]
\centering
\includegraphics[width=.8\textwidth]{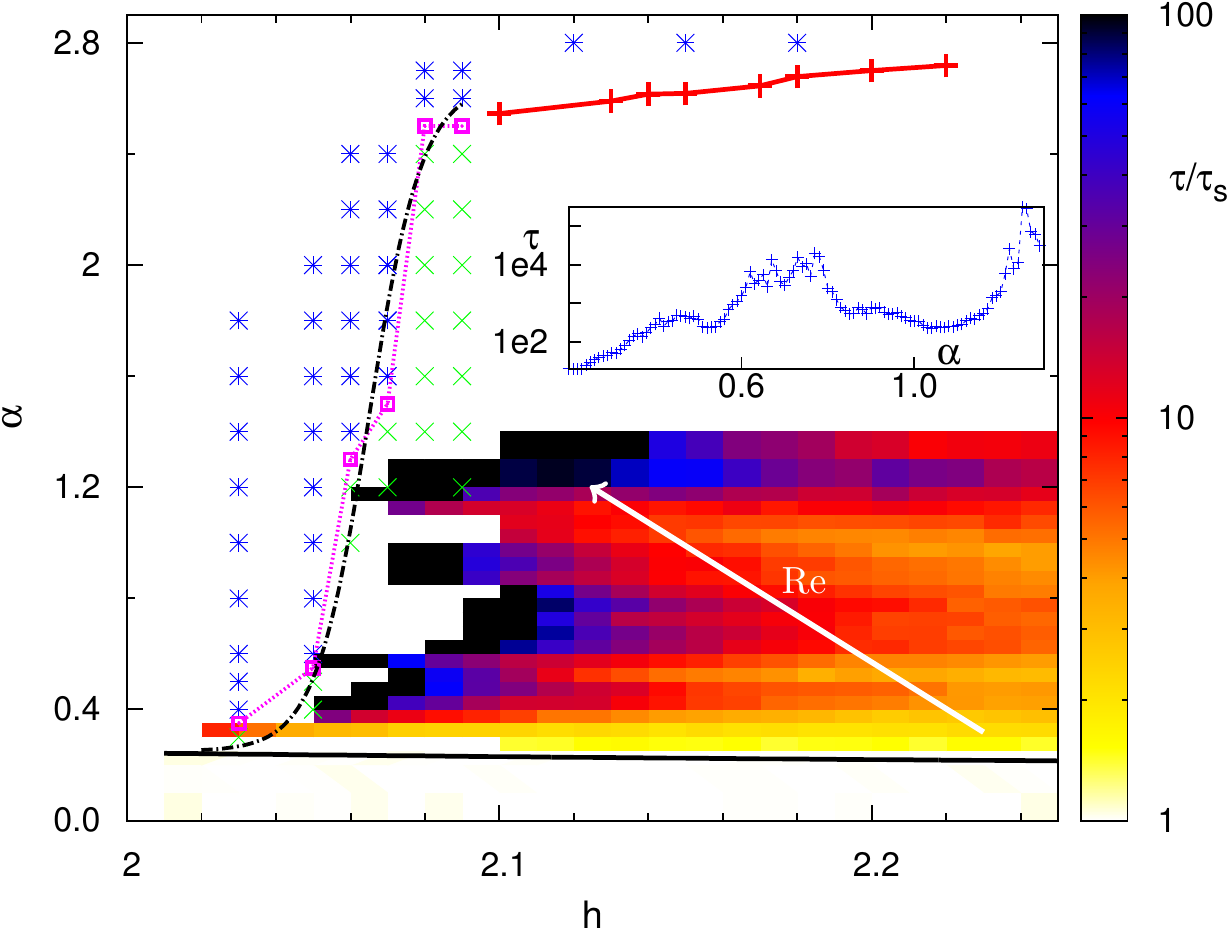}

\caption{Phase diagram for the qualitative behavior of chaotic
  transients of the UCML dynamics \Eq{model}. For
  $\alpha<\alpha_P\simeq 0.23$ (solid black line) chaotic behavior
  does not propagate along the system and transients have an average
  lifetime $\tau_s$.
  The lifetime of puffs is shown in false color according to the
  legend to the right of the plot. Green crosses mark puffs with a
  lifetime $\tau > 10^3 \tau_s$ from single trajectories, and slugs are indicated by blue
  stars.
  The dotted black line marks the transition from puffs to slugs upon
  decreasing $h$, and the solid red line is the border 
  $\alpha_{PS}$ as determined from \Fig{slug_velocities} (\cf main
  text). The white arrow indicates the trends how $\alpha$ and $h$ change upon increasing Reynolds number.
  The inset shows the lifetime $\tau$ as a function of $\alpha$ for
  fixed $h=2.1$.
  \label{fig:phase_space}}
\end{figure}

%
The key element is the possibility of a bifurcation analysis at an unphysical
  parameter value $\alpha_{sn}$ in order to gain insight into the
  physical transition from puffs to slugs.
By definition a slug amounts to a region in parameter space where the
leading edge moves faster than the trailing edge, $v_l - v_t > 0$ (see
\cite{DuguetWillisKerswell2010} for this definition). This
region is bounded by those combinations of parameters where both
velocities match. 
We note that 
 the leading edge velocity $v_l(\alpha)$ is a function of
$\alpha$ only, and the trailing edge velocity $v_t(h)$ depends solely
on $h$. 
Consequently, the transition line $\alpha_{PS}(h)$ 
can be found as root of the implicit function $0 = F(\alpha,h) = v_l(\alpha) - v_t(h)$,
and in view of  \Fig{slug_velocities} it can be parameterized by the velocity $v$.
We use a graphical algorithm based on \Fig{slug_velocities} to
  solve $v_l(\alpha) = v_t(h)$ for the implicitly-defined $\alpha$.
Selecting a velocity $v$ on the ordinate axis one can look up the
respective intersections with the data $v = v_l(\alpha)$ and
$v = v_t(h)$ and read off the corresponding values for
  $\alpha$ and $h$. Hereby, we obtain a function $\alpha(h)$.  

An example is given by the dotted gray line in \Fig{slug_velocities}
which starts off $v=0.18$, and then identifies the corresponding
control-parameter values of $v_l(\alpha \simeq 2.61) = 0.18$ and
$v_t(h \simeq 2.16) = 0.18$ for the desired leading and trailing edge
velocity, respectively.
This provides the point $\alpha_{PS}(h=2.16) \simeq 2.61$ of the
boundary $\alpha_{PS}(h)$ separating regions where puffs
and slugs are expected to arise. 
Varying $v$ results in an implicit description of the boundary
$\alpha_{PS}(h)$. 
The rule provides an accurate prediction for the transition
threshold where slugs arise when increasing $\alpha$ for a fixed value
$h \gtrsim 2.1$  (solid red line in \Fig{phase_space}). 
\subsection{Puffs}
\label{sec:puffs}
\begin{figure}[t]

  \centering \subfigure[Puff liftime vs.~single-site lifetime for
  \newline$\alpha < \alpha_p = 0.23$.]{\label{fig:lifetime_decay}\includegraphics[width=0.45\textwidth]{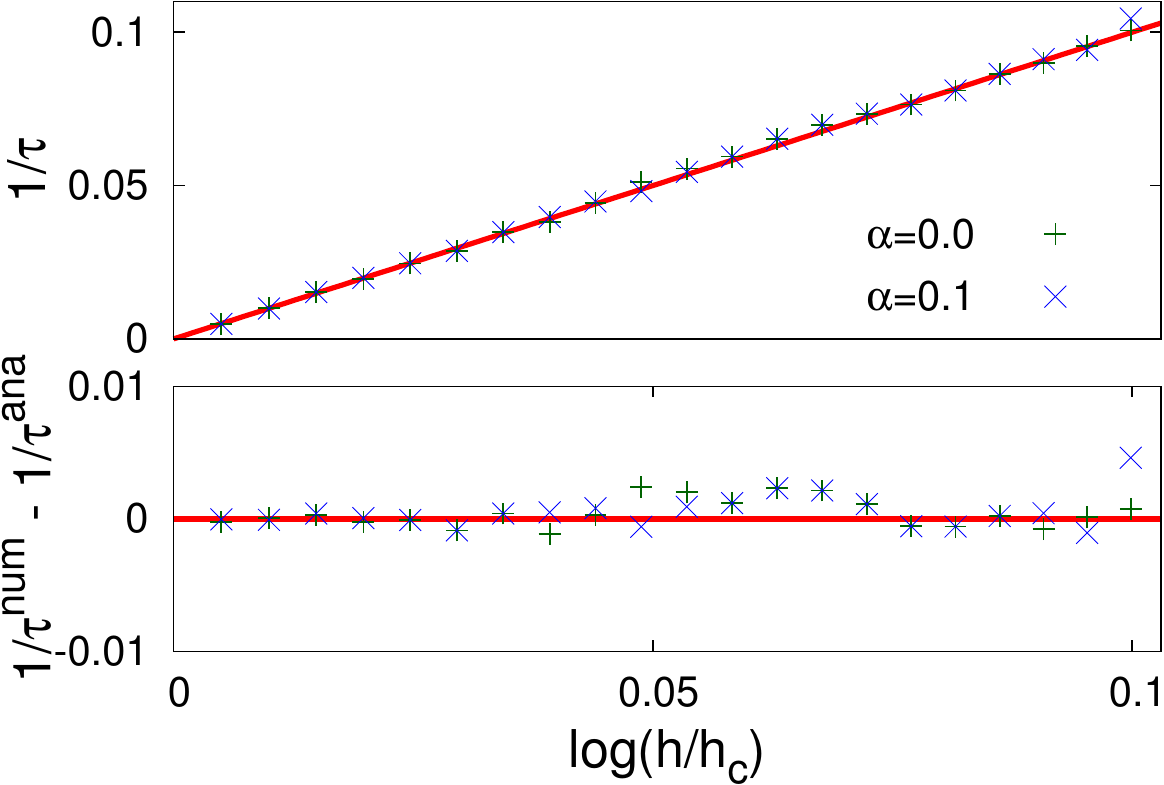}}%
  \hfill
  \subfigure[Critical curve for laminar-puff transition in parameter
  plane.]{\label{fig:alphacr}\includegraphics[width=0.45\textwidth]{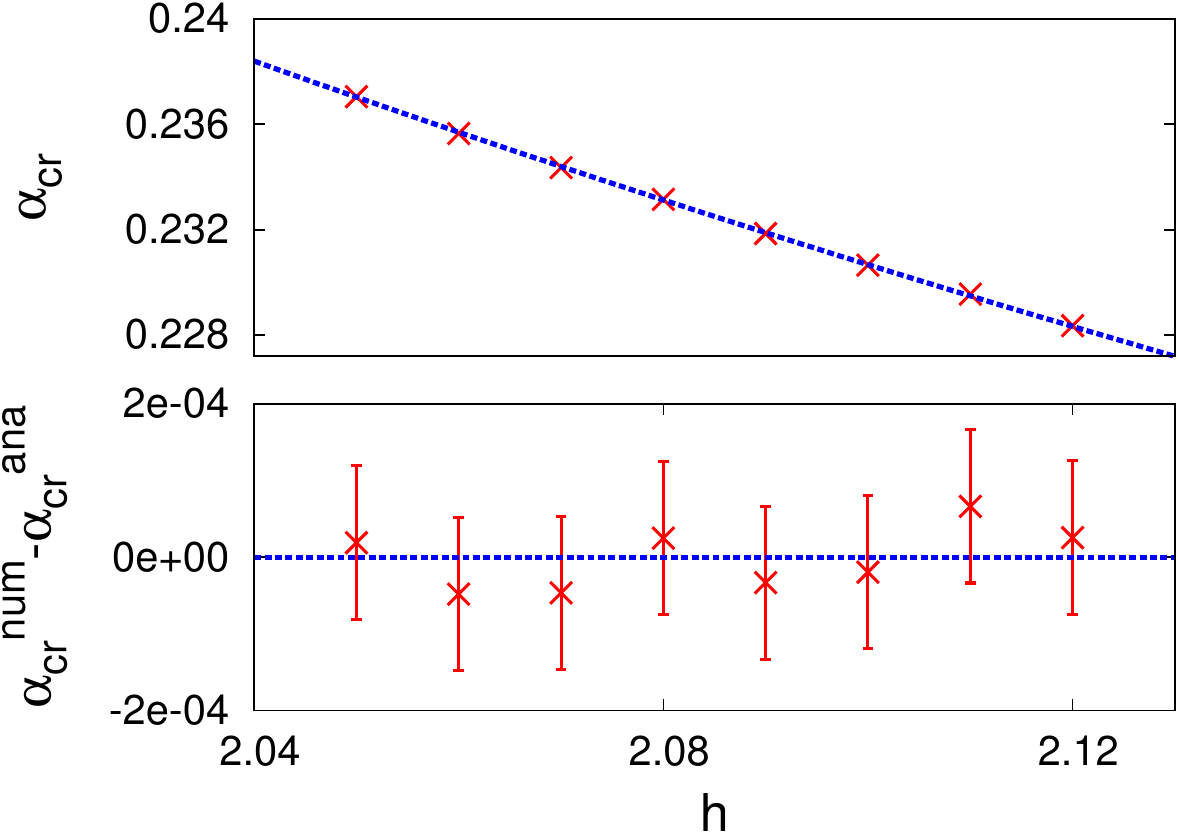}}
\caption{\textup{(a)}  For $\alpha < \alpha_{P}$ the lifetime follows exactly the
  prediction $\tau = \tau_s \equiv [\ln (h/h_c)]^{-1}$ provided by the
  solid red line. The relative difference to the analytical prediction
  is shown in the lower panel. \textup{(b)}
  Critical curve for the transition from a laminar state to a puff
  state. The implicitely-defined curve in \Eq{onset} is shown as
a blue dotted line. The results from numerical simulations (red crosses) are
in perfect agreement with the analytical result. To show the accuracy
of the prediction, the difference between the numerical results and analytical
prediction \Eq{onset} is shown in the lower panel. The errorbars
reflect the discretization of the numerical data set.} 
  \label{fig:lifetime_dec}
\end{figure}

Close to the boundary for the transition from puffs to slugs one observes puffs with an exceedingly small
probability to decay. Upon increasing $h$ and decreasing $\alpha$ the
average lifetime of puffs decreases as indicated by the color code in
\Fig{phase_space}. 
As long as $\alpha > \alpha_P \simeq 0.23$ (solid black line in \Fig{phase_space}) the average lifetime $\tau$
exceeds the one of the single-site dynamics by orders of magnitude;
its parameter dependence will be discussed in Section \ref{ssec:lifetime}.
Below the threshold one exactly recovers the single-site
lifetime as demonstrated in \Fig{lifetime_decay}.
The threshold coincides with the critical coupling strength where 
turbulence starts to propagate along the lattice, corresponding to the
emerging of the first unstable coherent states, \ie propagating
wave-like solutions of the Navier-Stokes equation, in pipe flow. It 
amounts to the solution of the implicit equation 
\begin{equation}
 \alpha_{P} \, g(x_2^f) = \delta
  \quad 
  \textrm{ involving the fixed point } 
  \quad
  x_2^f = \frac{h(2+\delta)}{h+1} 
  \label{eq:onset}
\end{equation}
where $f( x_2^f ) = x_2^f$. 
This equation expresses that the coupling reaches a level where repeated
kicks, that appear when $x^i$ takes the value of the right-most fixed
point of $f$, can push site $i+1$ from the laminar into the chaotic
state.
It is not possible to trigger
chaotic transients with the other fixed point $x_1^f$ of $f$, since
$g(x_1^f) <0$. 
The prediction \Eq{onset} is shown
by the solid black line marking the lower border of the colored region
in \Fig{phase_space}. This exact analytical result does not have any
free fit parameters. It exactly agrees with the numerically
determined threshold to an accuracy limited only by the numerical
error of the order of $10^{-4}$ (cf.\ \Fig{alphacr}).

\subsection{Lifetime of Turbulent Puffs}
\label{ssec:lifetime}
\begin{figure}[t]
  \centering 

\includegraphics[width=0.85\textwidth]{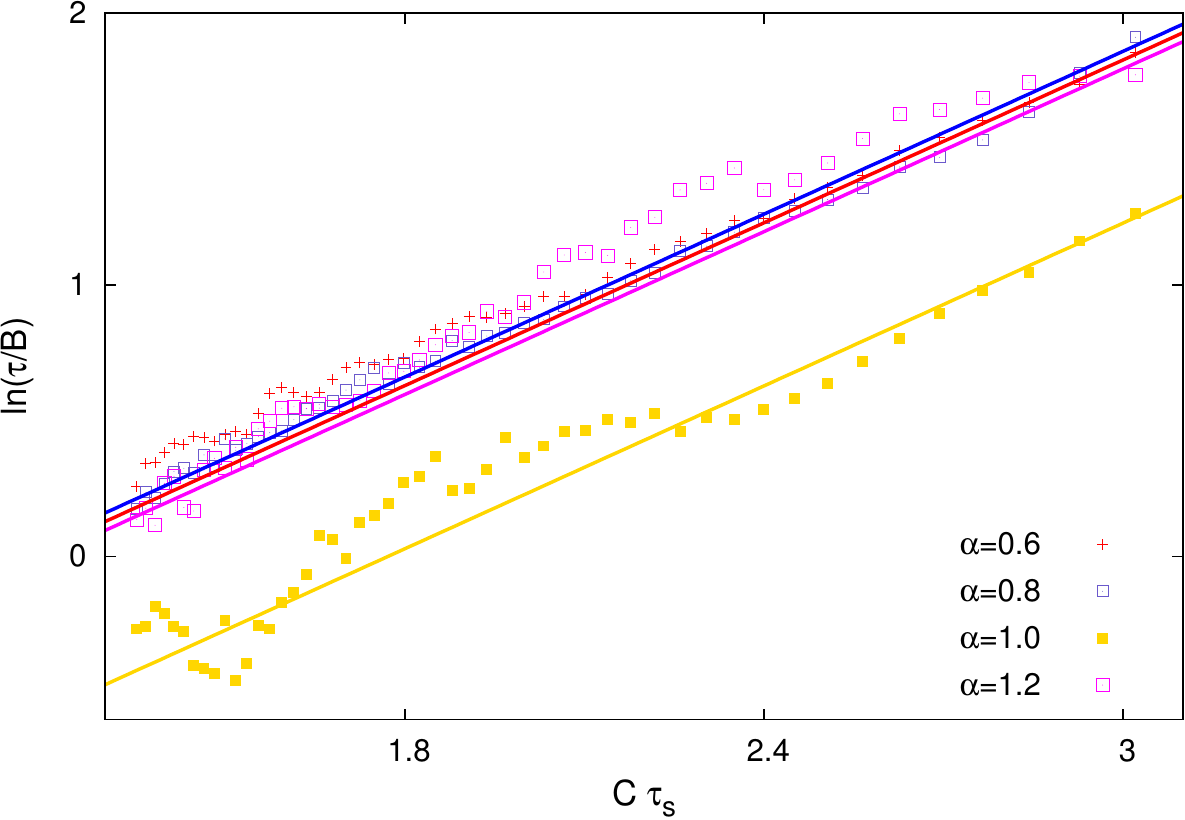}%
\caption{Increase of the average lifetime $\tau$ as function of 
  $\tau_s$. %
  Numerical data for different $\alpha > \alpha_{P}$ (different
  points as indicated in the legend to the lower right) agree very
  well with the super-exponential scaling, Eq.~(\ref{eq:superexp}).} 
\label{fig:lifetime}
\end{figure}
We now turn to the $h$ dependence of the lifetime $\tau$. 
For $\alpha>\alpha_{P}$ one observes a sharp increase of $\tau$, 
\begin{equation}
\ln \frac{\tau}{B} = C \tau_s,
  \qquad 
  \tau_s = \left[ \ln \frac{h}{h_c} \right]^{-1}
\label{eq:superexp}
\end{equation}
when $h$ approaches $h_c = 2$ from above (see \Fig{lifetime}). 
The exponential increase identifies the finite-size traveling
structures in this coupled map lattice as a new type of super-long
transients \cite{TelLai2008}.
As demonstrated in \Fig{lifetime} the dependence
\Eq{superexp} yields straight lines when $\ln(\tau/B)$ is plotted as
function of $C\tau_s$ in agreement with the numerical data.

\section{Discussion}
\label{sec:discussion}

We have introduced a coupled map lattice, \Eq{model}, with a
non-trivial unidirectional coupling, \Eq{g(x)}, between neighboring
lattice sites. This dynamics has several interesting solutions, in
particular in its transient behavior: (i) immediately decaying; (ii)
structures with a fixed width that travel forward through the lattice
(\Fig{puff} and \Fig{puffshort}); (iii) chaotic structures whose width
grows linearly in time (\Fig{slug}).  We provided an analytical
prediction of the lifetime of solutions of type (ii) which is in
perfect agreement with numerical simulations (\Fig{lifetime_dec} and
\Fig{lifetime}).  Furthermore, we presented a mechanism for the
transition between solutions of type (ii) and (iii), and we identified
it as an unbinding transition. All results were obtained analytically
as well as numerically. The analytical treatment revealed the
fundamental mechanisms, \ie the importance of the non-linear
  coupling, leading to different types of dynamics.

\subsection{Lifetimes of Puffs}

A simple statistical argument suggests that \Eq{superexp}
  might faithfully describe the lifetime of turbulent puffs. A
  dynamical systems representation of the transient turbulent dynamics
  should have $\Re^{9/4}$ active degrees of freedom which are sparsely
  coupled \cite{Grossmann2000}.
  Moreover, the escape rate $\kappa = \tau^{-1}$ amounts to the
  probability to enter the absorbing state where the motion proceeds
  directly to the laminar state \cite{Ott2002}. Due to the nonlinear
  coupling between the modes this involves a constraint on a
  noticeably portion of the degrees of freedom. Applying this scaling
  to the slices of pipe flow that are described by single cells of our
  UCML provides a scaling
\begin{equation}
  \tau^{-1} \simeq B' \: \exp\left[- C' \; \Re^\beta \right]
    \quad \text{with } 1 \lesssim \beta \leq 9/4 \,.
\label{eq:ourScaling}
\end{equation}
and appropriately adapted constants $B'$ and $C'$.

Remarkably, this dependence is compatible with the findings
of experimental pipe-flow data \cite{HofLozarKuikWesterweel2008} where
the authors supposed that their data are best fitted by 
\[
\tau^{-1} \simeq \exp\left[ -\exp( c_1 \Re + c_2 )\right]
\quad \text{with } c_1 = 0.0057 \text{ and } c_2 = 8.7 \,
\] 
while other adequate fits can also be obtained by other
superexponential functions, in particular, by 
\[
\tau^{-1} \simeq \exp\left[-(\Re/c)^n\right] 
\quad \text{with } c = 1549 \text{ and } n = 9.95 \,.
\]
It is straightforward to check that the expression provided in
\Eq{ourScaling} provides another very good fit of the experimental and
numerical data in the relevant range, $1650 < \Re < 2050$, of the
Reynolds numbers considered in \cite{HofLozarKuikWesterweel2008}.

\subsection{Comparison to Other Approaches Addressing Pipe Flow with Dynamical Systems
  Tools}

We now compare the
  main features obtained for the UCML with those of previously
  published research addressing pipe flow with dynamical systems
  tools, in particular to
the models of Barkley \cite{Barkley2011}, Allhoff and Eckhardt
\cite{Allhoff2012} and Sipos and Goldenfeld \cite{Sipos2011}. 
Barkley's model \cite{Barkley2011} employs ideas from excitable media. It comprises
a diffusively coupled map lattice in $(1+1)$ dimensions and two
variables with a non-linear coupling. The model is designed to mimic
main
features of pipe flow: turbulent puffs, puff splitting and
turbulent slugs. A comparison to direct numerical simulations of pipe flow 
reveals qualitatively the same results after fitting of the main
system parameter to experimental data. Unfortunately, the model has a
rather big number of parameters. They are not studied and
not compared to pipe flow, and severely hamper a theoretical analysis of the model.
The main advantage of the presented UCML model over the model
presented in \cite{Barkley2011} is the use of only a
  few parameters. This allowed us to provide
  an analytical treatment of its
transitions and the parameter dependence of the lifetimes of
  transient structures.
In principle we can also match the behaviour of pipe flow by
  relating the system parameters $\alpha$ and $h$ to $\Re$. The limit
  $h\to h_c = 2$ corresponds to a weaker decay to the laminar state,
  while an increase in $\alpha$ indicates a stronger non-linear
  coupling. To comply with experiments, one hence should
  simultaneously increase $\alpha$ and decrease $h$ to
  $h_c$. Moreover, $\delta$ and the slope of $f$ and $g$ at
  $x=0$ can be adjusted to account for the change of the stability of
  the laminar motion in pipe flow.

The models presented in \cite{Sipos2011} and \cite{Allhoff2012} are
directed percolation models and therefore stochastic by definition. 
Also there models were fitted
to pipe flow data by appropriately adjusting their
parameters. 
Since the Navier-Stokes equation is
deterministic the comparison only applies to statistical
quantities and does not give insight into the mechanisms from a
dynamical systems perspective. The model of \cite{Allhoff2012} is a $(1+1)$-dimensional
model with two parameters and the one in \cite{Sipos2011} is
$(3+1)$-dimensional and reduced to a $(1+1)$-dimensional model with
one parameter. While \cite{Allhoff2012} makes no prediction about the
relation of the two system parameters to $\Re$, \cite{Sipos2011} 
relates its control parameter $p$ to $\Re$ in order to discuss the
superexponential scaling of puff lifetime $\tau(p)$. Consequently, the
superexponential scaling is not shown in \cite{Allhoff2012}, but the
results are compared with mean field approximations. Puff splitting is
only observed in \cite{Sipos2011} but not studied in detail.

\section{Conclusion and Outlook}
\label{sec:conclusion}

Coupled map lattices with a unidirectional
non-linear coupling capture salient features of turbulent pipe flow in
terms of two parameters $\alpha$ and $h$:
\\
\textbf{(i) } in reminiscence to the low Re regime of pipe flow all
perturbations immediately decay to the laminar state when the coupling
constant $\alpha$ is small;
\\
\textbf{(ii) } for larger values of $\alpha$ and sufficiently large
$h$ the model shows chaotic regions of a finite width that propagate
down the pipe.  These puff-like structures feature a lifetime
statistics in accordance with super-long transients in dynamical
systems;
\\
\textbf{(iii) } decreasing $h$ leads to an unbinding transition of the
leading and trailing edge of puffs, similar to the transition from
puffs to slugs;
\\
\textbf{(iv) } there is an intermittent propagation of the leading
edge of puffs and slugs (see the discussion in Section
\ref{ssec:laminarpuffs}). These jumps have also been observed in real
pipe flow. They comply with space-time plots of direct numerical
simulations displayed in
Refs.~\cite{DuguetWillisKerswell2010,MoxeyBarkley2010,Avila2011}.

The UCML model uses less unknown parameters than the model presented
in \cite{Barkley2011}, and a better characterization is obtained by
analytical results.  All features were analyzed by explicit
analytical calculations. In particular, we provided predictions for
the parameter dependence of the transition lines and the lifetime of
puffs.
An interesting aspect of discussing the parameter dependences of
  our present UCML is that all its paramters have a well-defined
  physical interpretation for the pipe flow: the slopes of $f$ and
  $g$ at $x=0$ characterizes the linear stability of the laminar
  state, $\delta$ describes the stability against finite-size
  perturbations, $h$ relates to the leading Lyapunov exponent of turbulent
  states and $\alpha$ is the strength of the nonlinear forward forcing
  in turbulent states.  Hence, our analytical predictions connect the
  Re dependence of these features of the dynamics to the lifetime of
  transient structures and to the critical Re numbers that separate
  regimes where one observes stable laminar pipe flow, puffs, and
  slugs.
Qualitatively, we expect that the limit $h\to h_c = 2$ corresponds to
a weaker decay to the laminar state, while an increase in $\alpha$
indicates a stronger non-linear coupling. To comply with experiments,
an increase in \Re\ will therefore lead to a simultaneous increase of
$\alpha$ and decrease of $h$ towards $h_c$.

Although the choice of a unidirectional, nonlinear coupling describes many
features of pipe flow, it fails to faithfully cope with certain aspects of the dynamics. 
Arguably, the most prominent of these features is
puff splitting, which is an important ingredient in
the transition from puffs to slugs in real pipe flow. 
The onset of
puff splitting is partly observed in \Fig{spacetime}, but much more
analysis and appropriate modifications of the model will be necessary to
obtain reliable
results.
After all, models with only nearest-neighbor coupling along the pipe can not be expected to capture the
long-range coupling via the pressure that
relaminarizes the leading puff unless it is separated by a sufficient
distance \cite{Avila2011}. 

Besides gaining insight into pipe flow, it is also of interest to
  study the dynamics of the unidirectionally coupled map lattice on
  its own. Beyond the realm of diffusive coupling not many coupled map lattice models have been studied
  analytically.

  As an outlook, we mention that a probabilistic pendant to coupled
  map lattices is formalized in the theory of directed percolation
  \cite{Sipos2011} when one interprets the transition from puffs to
  slugs as threshold in directed percolation.  The
  application of percolation to advectively-diffusively coupled map
  lattices has been studied in
  \cite{ChateManneville1988,GrassbergerSchreiber1991,Ginelli2003}.  A
  comparison between our UCML and directed percolation invites further
  research. This future research can go into the direction of a
  systematic coarse graining to obtain the equivalent stochastic
  dynamics by means of symbolic dynamics. This coarse-graining
    would then provide a rigorous
  connection to percolation theory.

\section*{Acknowledgments}

  The authors acknowledge very inspiring discussions with Reiner Kree,
  Arkadii Pikovsky, Antonio Politi, Lamberto Rondoni and Tam\'as
  T\'el.  They are grateful to Kerstin Avila and Bj\"orn Hof for
  sharing data prior to publication, and to Lamberto Rondoni and
  Tam\'as T\'el for comments on a draft of the manuscript.


 \bibliographystyle{siam}
\bibliography{cml_turbulence}

\end{document}